\documentclass[a4paper, 12pt]{article}
\pdfoutput=1
\usepackage{amssymb}
\usepackage[utf8]{inputenc}
\usepackage{amsmath}
\usepackage{amsthm}
\usepackage{amstext}
\usepackage[brazil, english]{babel}
\usepackage{graphicx}
\usepackage{fancyhdr,setspace}
\usepackage{lastpage, caption, booktabs, newfloat, float, array}
\usepackage{ifthen}
\usepackage[colorlinks=true,linkcolor=black, urlcolor=blue]{hyperref}
\usepackage[headsep=0.5cm, hmargin=1in,centering,twoside,bindingoffset=0cm, left=3cm, right=3cm, top=2.5cm, bottom=3cm,]{geometry}
\usepackage{titling}
\usepackage{balance}
\usepackage{algorithmic}
\usepackage{textcomp}
\usepackage{setspace}
\usepackage[T1]{fontenc}
\usepackage{yhmath}
\usepackage{psfrag}
\usepackage{framed}
\usepackage{enumerate}
\usepackage{natbib}
\usepackage{multirow}
\usepackage{booktabs,siunitx,makecell}
\usepackage[referable]{threeparttablex}
\usepackage[shortlabels]{enumitem}
\usepackage{lettrine}
\usepackage{comment}
\usepackage{cases}
\usepackage{titlecaps}
\usepackage{hyperref}
\usepackage{caption}
\usepackage[export]{adjustbox}
\usepackage{capt-of}
\usepackage{varwidth}
\DeclareFloatingEnvironment[
  fileext   = loq,
  listname  = List of Charts,
  name      = Chart,
  placement = htbp,
]{quadro}
\addto\captionsbrazil{
}
\thanksmarkseries{alph}
\fancypagestyle{firstpage}{%
\fancyhf{} 
\fancyhead[L]{\emph{ \scriptsize 
\noindent Cadernos do IME - S\'{e}rie Matem\'{a}tica}\\[2pt] 
\emph{ \scriptsize \textbf{e-ISSN:} 2236-2797 (online)}}
\fancyhead[R]{\scriptsize\emph{DOI: 10.12957/cadmat.2020.XXXXX}\\[4pt]
\emph{Recebido em 00 Mês 2020  $\mid$  Aceito em 00 Mês 2020}}
}
\begin{document}
\title{{\bfseries \large AN OPTIMIZATION APPROACH TO STUDY THE PHASE CHANGING
		BEHAVIOR OF MULTI-COMPONENT MIXTURES}\\[1ex]
  \itshape \normalsize ESTUDO DO COMPORTAMENTO NAS MUDANÇAS DE FASE EM MISTURAS MULTI-COMPONENTES, UMA ABORDAGEM DE OTIMIZAÇÃO  \\}
\date{}
\author{\sc
GUSTAVO E. O. CELIS\thanks{Universidade Federal do Rio de Janeiro -- Engenharia Mecânica, Rio de Janeiro - Brasil; ORCID: 0000-0002-4042-9093 \textbf{E-mail:} \href{mailto:gustavo.oviedo@ufrj.br}{gustavo.oviedo@ufrj.br}}
\and
REZA AREFIDAMGHANI\thanks{Universidade Federal do Rio de Janeiro -- Instituto de Matemática Aplicada, Rio de Janeiro - Brasil; ORCID:0000-0001-9885-2532 \textbf{E-mail:} \href{mailto:reza.arefi@matematica.ufrj.br}{reza.arefi@matematica.ufrj.br}}
\and
HAMIDREZA ANBARLOOEI\thanks{Universidade Federal do Rio de Janeiro -- Instituto de Matemática Aplicada, Rio de Janeiro - Brasil; ORCID:0000-0002-3220-7255 \textbf{E-mail:} \href{mailto:hamid@matematica.ufrj.br}{hamid@matematica.ufrj.br}}
\and
DANIEL O. A. CRUZ\thanks{Universidade Federal do Rio de Janeiro -- Engenharia Mecânica, Rio de Janeiro - Brasil; ORCID: 0000-0003-2465-5119 \textbf{E-mail:} \href{mailto:doac@mecanica.coppe.ufrj.br}{doac@mecanica.coppe.ufrj.br}}
}

\maketitle
\thispagestyle{firstpage}

\vspace*{-2cm}
\begin{otherlanguage}{brazil}
\begin{abstract}
\singlespacing 
O dimensionamento, construção e operação dos processos de captura e armazenamento de carbono (CCS) e de recuperação avançada de petróleo (EOR) requerem de um entendimento detalhado do comportamento das fases resultantes em misturas multi-componentes de hidrocarbonetos e dióxido de carbono ($CO_2$), isto, sob condições de reservatório.  Para modelar esse comportamento, um sistema não linear que consiste de uma equação de estados e algumas regras de mistura (para cada componente) precisam ser resolvidas simultaneamente.  A mistura geralmente requer modelar a interação binária entre os componentes da mistura. Este trabalho emprega técnicas de otimização para aprimorar as previsões desse modelo, otimizando os parâmetros de interação binária. Os resultados mostram que os parâmetros otimizados, embora obtidos matematicamente, estão em faixas físicas e podem reproduzir com sucesso as observações experimentais, especialmente para os sistemas de hidrocarbonetos multicomponentes contendo dióxido de carbono nas temperaturas e pressões do reservatório.\\

\noindent{\bf Palavras-chave:} Predictive Peng–Robinson, Modelagem,  Mudança de fases, Hidrocarbonetos, $CO_2$, Otimização.
\end{abstract}
\end{otherlanguage}

\begin{otherlanguage}{english} 
\begin{abstract}
\singlespacing 
The appropriate design, construction, and operation of carbon capture and storage (CCS) and enhanced oil recovery (EOR) processes require a deep understanding of the resulting phases behavior in hydrocarbons-$CO_2$ multi-component mixtures under reservoir conditions. To model this behavior a nonlinear system consists of the equation of states and some mixing rules (for each component) needed to be solved simultaneously. The mixing usually requires to model the binary interaction between the components of the mixture. This work employs optimization techniques to enhance the predictions of such model by optimizing the binary interaction parameters. The results show that the optimized parameters, although obtained mathematically, are in physical ranges and can reproduce successfully the experimental observations, specially for the multi-component hydrocarbons systems containing Carbon dioxide at reservoir temperatures and pressures.\\



\end{abstract}
\end{otherlanguage}

\section{Introduction}
Many modern chemical and mechanical process designs require models which are capable of
predicting the equilibrium properties of mixtures (e.g. oils with different components consisting of Methane ($CH_4$), Ethane ($C_2H_6$), Carbon dioxide ($CO_2$)). Usually the available experimental data is limited for such mixtures and the model also should produce accurate predictions for the sub-critical and critical conditions where performing experiments (if possible) is difficult. This type of modeling is challenging physically and numerically. Over the last 25 years many researchers have been working on different aspects of this problem. The governing equation of the problem is based on an Equation of States (EOS) for every phase present in the mixture and a mixing rule to explain the interaction between these components. During many years, the classical approaching of van der Waals mixing rules used the binary interaction parameter ($k_{ij}$) as an input obtained from the experiments. Nowdays this approach  has loosed importance and new theoretical models  for the prediction of $k_{ij}$ has been proposed for the description of vapor-liquid equilibrium (VLE)  using cubic EOS. More information
and comparison between all these methods can be found in ([\textbf{1, 2, 3, 4, 5}] and [\textbf{8}]). Mathematically speaking, this type of modeling results in a highly non-linear system of algebraic equations. To find the physical solution of this system, one needs to use sophisticated iterative methods which depend strongly on the initial guess (providing proper initial guess which is also physical is an open question). The mathematical aspects of numerical solution method can be find in ( [\textbf{5, 6}] and [\textbf{9}]).

In this research we have investigated the behavior of the hydrocarbon mixtures with more than 20
components. In the present work the Peng-Robinson78 (PR78) equation of state with the group contribution formulation for the binary interaction parameters are to model the interaction between the components. The binary
interaction parameters (in the order of $n^2$, where $n$ is the number of components) play an
important role in the modeling and express the effect of each component on the others.
Unfortunately, the physical theory behind these parameters is not adequate and one needs to adjust
them using other knowledge (experimental data). 

Comparing the results of the modeling with available experimental data shows
discrepancies when $CO_2$ is present in the mixture. Using some physical intuition, it was realized before that adjusting the binary interaction parameter between $CO_2$ and $CH_4$ improves the results [\textbf{1}]. However, by adjusting this single parameter at each temperature, a slight agreement was observed at the experimental bubble points, while quite big discrepancies remained at the dew points. The main idea of the present work is to adjust the binary interaction parameters to minimize this difference. This work employs optimization techniques to optimize the predictions of the phase behavior of systems containing $CO_2$ and multi-component hydrocarbons at reservoir temperatures and pressures. \\

The objectives of this work are divided into three main parts: 
First, to study the phase changing behavior of ($CO_2$ + synthetic crude oil) over wide ranges
of temperature, pressure and Gas-Oil Ratio (GOR) by introducing a successive method for starting the phase equilibrium calculation algorithm with a suitable initial guess (PPR78+InG). This modification led into improvements of our proposed method (PPR78+InG) by reducing the cost function MSE (Mean Squared Error) compared to the  other results reported in literature  [\textbf{1}]; Second, to adapt an optimization method to optimize the binary interaction parameters (interaction between $CO_2$ and other components) and optimize the cost function (PPR78+Opt); and third,to compare the outcomes of our suggested approaches (PPR78+InG and PPR78+Opt) with those of other known predictive models, namely, predictive Peng Robinson and a modified version of it with an adjustment strategy  [\textbf{1}], denoted by PPR78 and PPR78+Adj respectively.

\section{Governing Equations}
The PPR78  [\textbf{3}] cubic EoS is an extension of the PR78.  For a single component and using the so-called Soave alpha function, \(\alpha_i(T)\), the PR78 model is given by  the Eq. \ref{eq:eos} 
\begin{equation}\small
		P = \frac{RT}{v -b_i} - \frac{a_i(T)}{v(v+b_i)+b_i(v-b_i)}
	\label{eq:eos}
\end{equation}
{\small
	\[\left\{
	\begin{array}{lr}
		R = 8.314472 \left[\frac{\mathrm{J}}{\mathrm{mol}.\mathrm{K}}\right], \quad   \qquad \qquad \qquad \qquad \qquad \qquad \qquad a_i=0.457235529 \displaystyle\frac{R^2T^2_{c,i}}{P_{c,i}}\alpha_i(T), \\
		b_i= 0.0777960739 \displaystyle\frac{RT_{c,i}}{P_{c,i}}, \quad  \qquad \qquad \qquad \qquad   \alpha_i(T)=\displaystyle\frac{R^2T^2_{c,i}}{P_{c,i}}\left[1+m_i \left(1-\sqrt{\frac{T}{T_{c,i}}}\right)\right]^2, \\
		\\
		m_i = 0.37464+1.54226w_i-0.26992w^2_i  \qquad \qquad \qquad \qquad \qquad \qquad  \qquad \quad \textrm{if} \ w_i\leq 0.491\\
		m_i = 0.374642+1.48503w_i-0.164423w^2_i + 0.016666w^3_i \qquad \qquad \qquad \qquad \textrm{if} \ w_i> 0.491 
	\end{array}
	\right.
	\]}
where $P$ is the pressure, $R$ the ideal gas constant, $T$ the temperature, $a$ and $b$ are the energy and co-volume parameters respectively.  $v$ is the molar volume, $T_c$ and $P_c$ the critical temperature and pressure respectively, and $w$ the acentric factors.  To apply such an EoS to the mixtures, mixing rules need to be used to calculate the $a$ and $b$ values related to the mixture.  The classical van der Walls mixing rules approach  used in this study are:
		 \begin{equation}
	 	\begin{split}
		 		a&= \displaystyle\sum_{i=1}^{N}\sum_{j=1}^{N} z_i z_j \sqrt{a_i a_j} (1-k_{ij}(T)) \\
		 		b&= \displaystyle\sum_{i=1}^{N} z_i b_i,
		 	\end{split}
	 	\label{eq:ab}
	 \end{equation}
where $z_i$ and $z_j$ represents the molar fraction of every component in a mixture, $N$ is the number of components in the mixture, and  \(k_{ij}(T)\)  is the binary interaction parameter between molecules \(i\) and \(j\).  When \(i=j\), \(k_{ij}=0\), and also \(k_{ij}=k_{ji}\).  Usually \(k_{ij}\) is fitted experimentally, and it has been shown that it is temperature dependent.  Based on the group contribution method (GCM), \(k_{ij}(T)\) to the PPR78 model is given by
\begin{equation}\small
	\begin{split}
		k_{ij}(T)=\frac{-\displaystyle\frac{1}{2}\displaystyle\sum\limits_{i=1}^{N_g}\sum\limits_{j=1}^{N_g} (\alpha_{ik}-\alpha_{jk})(\alpha_{il}-\alpha_{jl})A_{kl}. \left(\frac{298.15}{T}\right)^{\left(\frac{B_{kl}}{A_{kl}}-1\right)}-  \left(\frac{\sqrt{a_i(T)}}{b_i}-\frac{\sqrt{a_j(T)}}{b_j}\right)}{2\displaystyle\frac{\sqrt{a_i(T).a_j(T)}}{b_i.b_j}},
	\end{split}
	\label{eq:kij}
\end{equation}
where \(a_i\) and \(b_i\) are values obtained from mixing rules, and \(N_g\) is defined as the number of different groups present in the molecule.  $\alpha_{ik}$, $A_{kl}$ and $B_{kl}$ are associated parameters related to the GCM; for more information, we refer the reader to   [\textbf{5}].

To fully explain the system, one needs to find the phases present in the mixture. This can be done using Flash calculation. The method is primarily concerned with phase molar fraction ($\beta_i$) determination, equilibrium phases number and the molar composition of every component (\(x_i\)) in each phase (\(X_k=\sum{x_i}\)) 
\citep{Michelsen2007}.  For a given global composition of a mixture (\(Z=\sum{z_i}\)) and two thermodynamic independent intensive properties,  often pressure (\(P\)) and temperature (\(T\)) is possible to perform the Flash calculation.  In ($PT$)-Flash, the solution corresponds to the global minimum in Gibbs Energy of the multiphase system as necessary and sufficient condition for equilibrium.  The ($PT$)-Flash
methodology selected in the present work was the one proposed by  [\textbf{9}].

The above model depends on $k_{ij}(T)$, the binary interaction parameter between components $i$ and $j$. Eq. \ref{eq:kij} provides an estimate for this parameter. However, it is not accurate enough, specially in the presence of $CO_2$ in the mixture. In the next section, we use optimization technique to find an optimized value for $k_{ij}$ around the value predicted by Eq. \ref{eq:kij}. Then, we show that these optimized values reproduce better the experimental observations.

\section{ EVALUATION METRIC}\label{Metrics}
Evaluation metric refers to a measure that we use to evaluate different models. Prior to any modeling process, selecting an appropriate evaluation metric is a crucial decision that demands for a complete understanding of the project's goals.
A cost function is used to gauge the performance of the proposed model (Machine Learning model or predictive model). A predictive model without a cost function is  pointless, since the cost function is useful for determining how well a proposed model performs. The appropriate selection of the cost function contributes to the model's credibility and reliability.
There are a number of popular cost functions, namely Mean Absolute Error (MAE), Mean Squared Error (MSE), Root Mean Squared Error (RMSE), Root Mean Squared Logarithmic Error (RMSLE), and Maximum Error (MAXE).

To begin, consider the Mean Absolute Error (MAE), which is the mean absolute difference between the actual and predicted values. It is expressed by $     \text{MAE}=\frac{1}{n} \left[\sum\limits_{i=1}^n |p_{_i}-{a_i}| \right]$,
where $p_i$ is the predicted value of instance $i$, $a_i$ is real (target) value of instance $i$ and $n$ is the total number of instances.
MAE is resistant with respect to outliers, because it does not punish high mistakes brought on by outliers, meaning  it is insensitive to outliers.

The maximum error function (MAXE) calculates the highest residual error. A measure that captures the worst-case difference between the predicted and actual values.
The Mean Squared Error (MSE), which measures the difference between desired and predicted solutions, is one of the most popular metric. The optimized cost function can only be obtained with the smallest MSE value. The MSE is defined as, $ \text{MSE}=\frac{1}{n} \left[\sum\limits_{i=1}^n (p_{i}-{a_i})^2\right]$.    

The main draw for using MSE is that it squares the error, which results in large errors being punished or clearly highlighted. It’s therefore useful when working on models where occasional large errors must be minimised.
The fact that the order of loss is greater than that of the data is another problem with MSE. The data being of order one and the cost function (MSE) has an order of two. Therefore, we are not able to immediately correlate the data to the error. Hence, we take the root of Mean Squared Error: $\text{RMSE}=\sqrt{\frac{1}{n} \left[\sum\limits_{i=1}^n (p_{i}-{a_i})^2\right]}$.

Note that we are not modifying the cost function in this case, and the solution remains the same. All we have done is reduce the order of the cost function by taking the root. 
RMSE can be employed when we want to penalize high errors; in other words, The square root in RMSE makes sure that the error term is penalized. 
Note that both RMSE and MSE measure relative errors with the same intensity, since the square root in RMSE is taken over the entire sum. The use of RMSE is justified because it has the same unit as the variables.

The next metric is known as Root Mean Squared Logarithmic Error (RMSLE), and it is quite similar to RMSE except that the log is applied before computing the difference between actual and predicted values. Both big and small errors are treated similarly.  Compared to RMSE, RMSLE is less sensitive to outliers. In the case of RMSE, the presence of outliers can explode the error term to a very high value. But, in the case of RMLSE the outliers are drastically scaled down therefore nullifying their effect. Thus, it relaxes the penalty of high errors resulting from the presence of the log. It is formulated as: $    \text{RMSLE}=\sqrt{\frac{1}{n} \left[\sum\limits_{i=1}^n \left(\log(1+p_{i}))-\log ({1+a_i})\right)^2\right]}$.

\section{RESULTS}

In this section, various metrics are employed to compare the performance of of the previously described methodologies, namely, PPR78, PPR78+Adj  [\textbf{1}] and PPR78+InG, PPR78+Opt (present work). To begin with, the experimental result reported in  [\textbf{1}] is utilized as a benchmark to examine different approaches in predicting these results.

The problem under consideration is an classic constrained optimization defined as: minimizing the cost function on the condition of binary interaction of $CO_2$-$CH_4$, being a number positive or negative near zero. The grid search optimization through a simple implementation of it in Python is utilized to solve all of the problems in terms of different evaluation metrics and temperatures. First, a grid search is employed to minimize the cost function over the interval [-0.2, 0.2]. The $k_{CO_2-CH_4}$(with the precision 0.01) associated with the smallest cost function is selected as the candidate of the neighborhood of the optimal solution. Next, another grid search is performed over the interval [$(k_{CO_2-CH_4}-\epsilon$, $k_{CO_2-CH_4}+\epsilon$], where $\epsilon=0.01,$ is applied to find the optimal value for $k_{CO_2-CH_4}$ with the precision $0.001,$ i.e., the minimizer of the cost function.

The cost function values MAE, MSE, RMSE, RMSLE, and MAXE produced by different predictive models for Live Oil 1 and Live Oil 2 (introduced in  [\textbf{1}])  and also for various temperatures 323.15 K, 373.15 K and 423.15 K   are provided separately in Table \ref{table1 :costFunctions}.  These evaluation metrics (as described in section \ref{Metrics}) are utilized for each approach, and the
best performance emphasized in bold. However, as indicated before, our focus will be on the
three measurements MSE , RMSE and RMSLE. Nevertheless, as can be seen in in Table \ref{table1 :costFunctions}, according to all metrics, our suggested model (PPR78+Opt) is significantly better than other approaches. That is, the cost function derived by PPR78+Opt for Live Oil 1, Live Oil 2, and varied temperatures is much smaller than that obtained by other methods.

\begin{table}[h]\tiny
	\renewcommand{\arraystretch}{1.5}
	\centering
	\caption{\small The cost function values provided by different predictive models, separately in various temperatures.}
	\label{table1 :costFunctions}
	\begin{tabular}{|l|l|rrrrr|rrrrr|}
		\hline 
		\textbf{Temp [K]} & \textbf{Method}  & \multicolumn{5}{c|}{\textbf{  Live Oil 1} }  & \multicolumn{5}{c|}{\textbf{Live Oil 2} }      \\ \hline
		\multirow{5}{*}{\textbf{323.15 K}} & \multirow{2}{*}{}         & MAE & MSE &  RMSE &  RMSLE & MAXE                  & MAE   & MSE    &  RMSE & RMSLE & MAXE   \\ \cline{3-12} 
		
		&  \textbf{PPR78}       & 4.489   & 23.376   & 4.834   & 0.229   & 6.949     & 5.448 & 31.408 & 5.604 & 0.270 & 7.188 \\
		& \textbf{PPR78+Adj}     & 2.649   & 11.211   & 3.348   & 0.167   & 6.297     & 3.247 & 14.560 & 3.815 & 0.204 & 5.922  \\
		& \textbf{PPR78+InG}   & 3.355   & 13.494   & 3.673   & 0.180   & 6.058     & 5.281 & 29.221  & 5.405 & 0.303 & 7.319   \\
		& \textbf{PPR78+Opt}    & \textbf{0.73}& \textbf{1.44}& \textbf{1.201}& \textbf{0.064}& \textbf{3.07}& \textbf{1.13}& \textbf{1.71}&\textbf{1.3}& \textbf{0.055} &  \textbf{2.59}  \\
		\hline  
		\multirow{5}{*}{\textbf{373.15 K}} & \multirow{2}{*}{}  & MAE & MSE &  RMSE &  RMSLE & MAXE & MAE & MSE &  RMSE &  RMSLE & MAXE   
		\\ \cline{3-12} 
		&   \textbf{PPR78}      & 3.094   & 10.883   & 3.299   & 0.130   & 5.094      &  3.726 & 14.455 & 3.802  & 0.143  & 4.975 \\
		&  \textbf{PPR78+Adj}    & 2.213   & 6.794    & 2.606   & 0.101   & 4.768      &  2.087 & 6.312  & 2.512  & 0.107  & 3.946  \\
		& \textbf{PPR78+InG}  & 3.099   & 10.978   & 3.313   & 0.130   & 4.391       &  3.603 & 13.657 & 3.695  & 0.139  & 4.773   \\
		& \textbf{PPR78+Opt}   & \textbf{0.59}    & \textbf{0.63}     & \textbf{0.79}    & \textbf{0.034}   & \textbf{1.69}       &   \textbf{0.91} & \textbf{0.96}   & \textbf{0.98}   & \textbf{0.034}  & \textbf{1.37}     \\        \hline
		\multirow{5}{*}{\textbf{423.15 K}} & \multirow{2}{*}{}            & MAE & MSE &  RMSE &  RMSLE & MAXE                   & MAE   & MSE    &  RMSE  &  RMSLE  & MAXE   \\ \cline{3-12} 
		&  \textbf{PPR78}      & 2.579   & 8.406  & 2.899  & 0.097  & 4.321          & 3.606 & 14.184 & 3.766  & 0.127   & 5.963   \\
		& \textbf{PPR78+Adj}    & 2.002   & 6.276  & 2.505  & 0.082  & 4.184          & 2.171 & 8.876  &  2.979 & 0.109   & 6.494   \\
		& \textbf{PPR78+InG}  & 2.637   & 0.733 & 2.955  & 0.098  & 4.324          & 3.528 & 13.994 & 3.740  & 0.127   & 6.099   \\
		& \textbf{PPR78+Opt}   & \textbf{0.54}    & \textbf{0.35}   & \textbf{0.59}   & \textbf{0.021}  & \textbf{1.12} & \textbf{0.547} & \textbf{0.38}   & \textbf{0.616}  & \textbf{0.018}   & \textbf{1.07}    \\ \cline{1-12} 
		
	\end{tabular}
\end{table}

\medskip

\begin{figure}[h]
	\centering  \includegraphics[width=0.9\textwidth,trim={3cm 1cm 3cm 2cm} ]{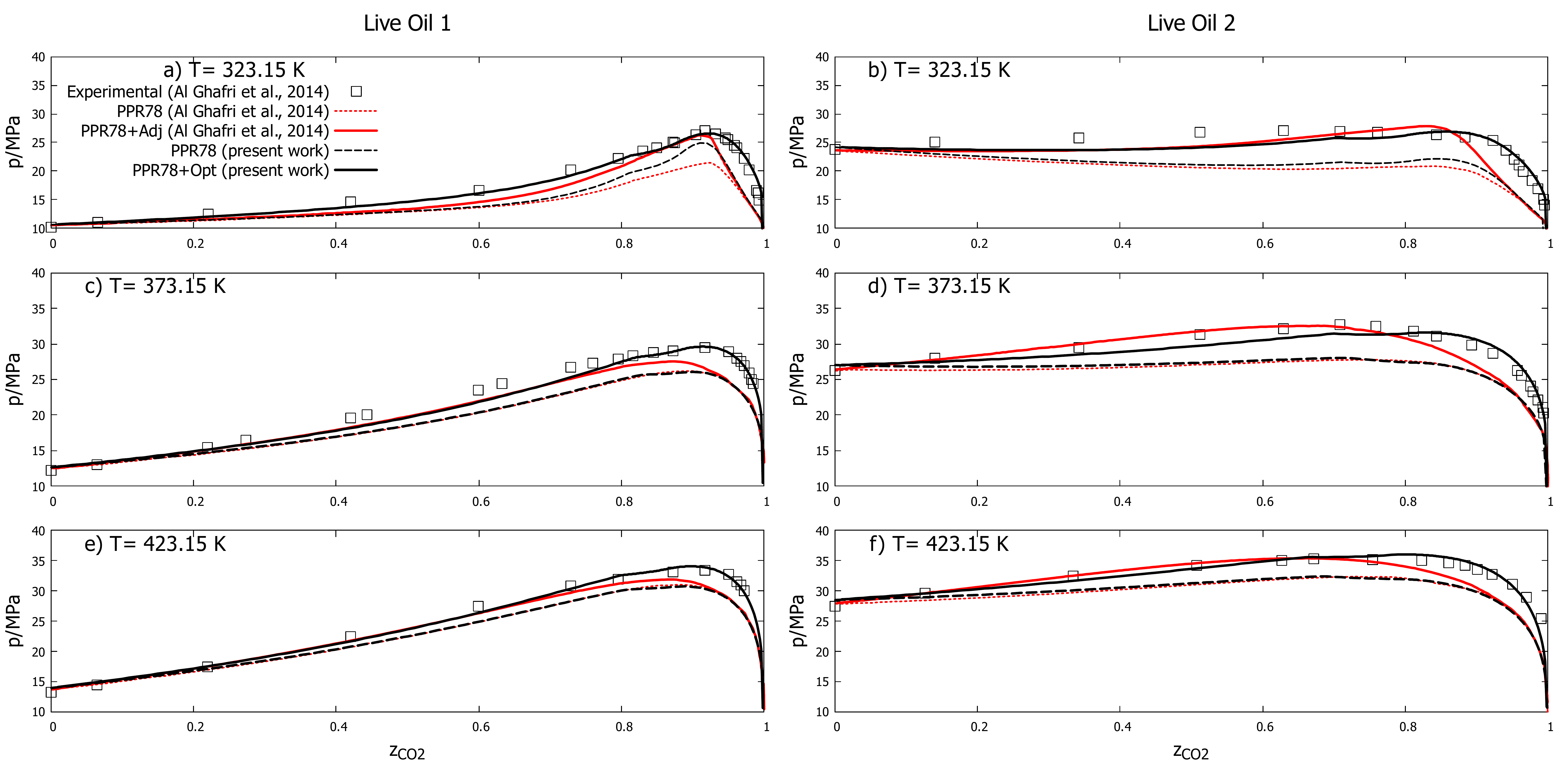}
	\caption{{\small Bubble-point and dew-point pressures ($p$) as a function of the mole fraction (z$_{CO_2}$ ) of
			$CO_2$ for ($CO_2$ + Live oil 1) figures (a, c, e), and ($CO_2$ + live oil 2) figures (b, d, e ):  T = 323.15 K;
			, T = 373.15 K and , T = 423.15 K.  }}
	\label{fig:6fig}
\end{figure}

The visualization of the results reported in Table \ref{table1 :costFunctions}, with respect to the prediction of binary interaction parameters for $CO_2$–$CH_4,$ can be seen in Figure \ref{fig:6fig}. As can be seen, the PPR78 predictive model has been improved by decreasing the cost function after being furnished with an initial guess (PPR78+InG), i.e., the black dashed lines are closer to the experimental data than the red dashed lines. There is a reason why this improvement is readily seen lower temperatures but as those aspects require more technical information, we omit it. Another important observation in Figure \ref{fig:6fig} is that predictive models PPR78, PPR78+Adj  [\textbf{1}] and PPR78+InG are accurate at predicting bubble points but less so at predicting dew points. However, our proposed optimization approach (PPR78+Opt) accurately predicts the experimental results in both the bubble and dew points. 
Note that the cost functions obtained by PPR78+Opt for various mixtures and temperatures (reported in Table \ref{table1 :costFunctions}) are reasonably small, resulting in agreement of black solid lines with experimental points in Figure \ref{fig:6fig}.
It should be noted that this observations contributes to the results shown in Table\ref{table1 :costFunctions}, meaning that PPR78+Opt provides the smallest cost function when compared to other predictive methods. This demonstrates how successfully the proposed optimization method predicts the binary interaction between $CO_2$ and $CH_4$ and consequently reducing the cost function.

\begin{figure}[h]
	\begin{minipage}[b]{0.5\linewidth}
		\centering

			\captionof{figure}{\small The mean of {\textcolor{blue}{MAE}}, {\textcolor{magenta}{MSE}}, {\textcolor{red}{RMSE}}, {\textcolor{red}{RMSLE}}, {\textcolor{black}{MAXE}}, of all the methods, reported in Table  \ref{table 2:Maen Of Metrics}.}
			\label{fig:spiderweb}
			\captionof{table}{\small The mean of the metrics values of all employed methods}
			\footnotesize\centering
			\begin{tabular}{|l|c|c|c|c|c|}
				\hline
				Method     & MAE         & MSE             & RMSE           & RMSLE      & MAXE           \\\hline
				PPR78      & { 3.823} & { 17.118}    & { 4.034}    &   0.166    & { 5.748}     \\
				PPR78+Adj    & { 2.394} & 9.004           & 2.960         & 0.128      &   5.268               \\
				PPR78+InG  & { 3.583} & 13.679           & 3.796          & 0.162     &    5.198               \\
				PPR78+Opt    & {\bf 0.741 }     &  {\bf0.911}          &{\bf  0.912}         &{\bf  0.037}    &  {\bf    1.818  }          \\\hline 
			\end{tabular}
			\label{table 2:Maen Of Metrics}
	\end{minipage}\hfill
	\begin{minipage}[b]{0.4\linewidth}
		\centering
		\includegraphics[scale= 0.4, trim={10cm 6cm 8cm 7cm}]{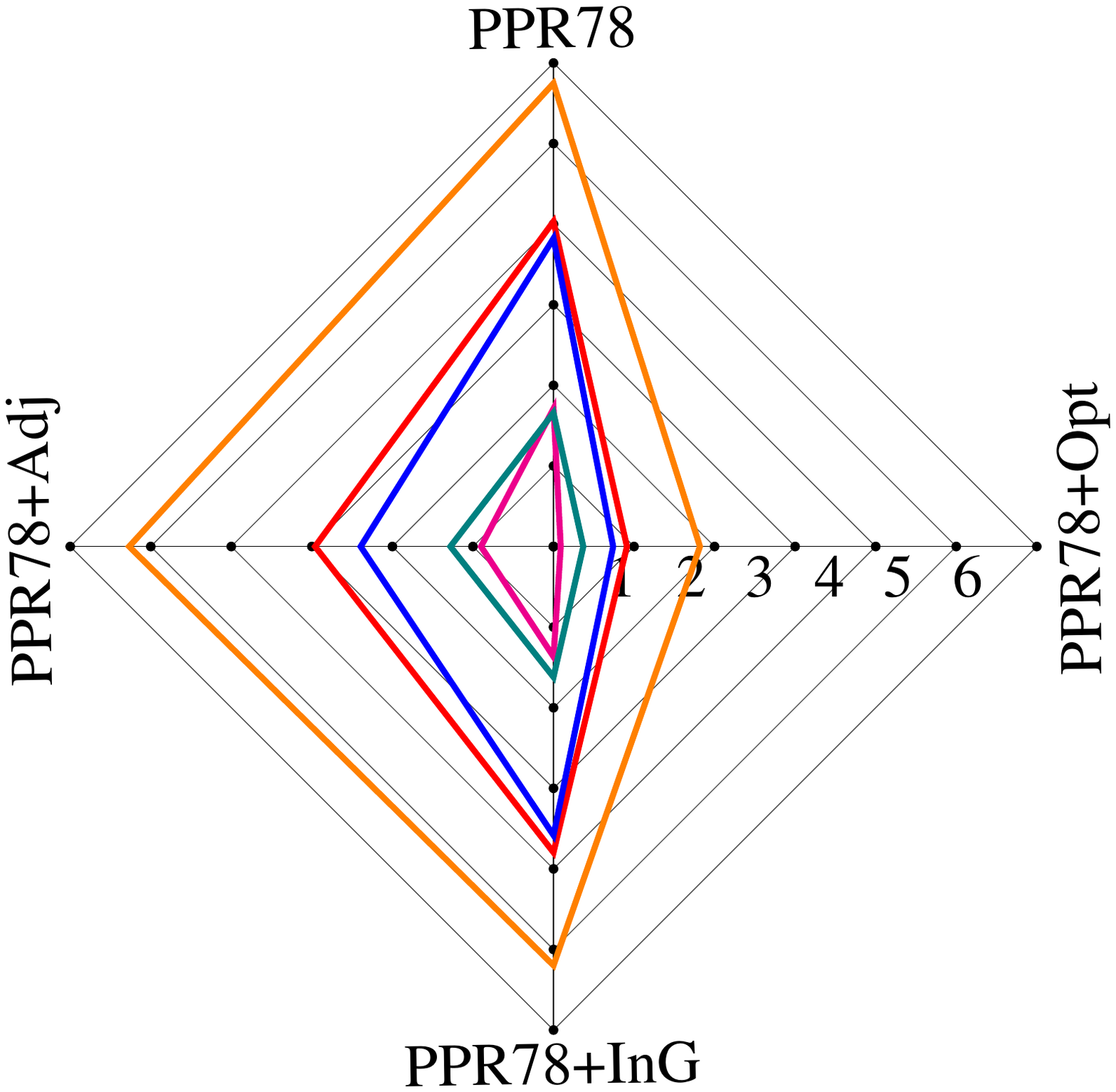}
	\end{minipage}
\end{figure}

To have another trustworthy measurement of the predictive approaches, the mean  values of MAE, MSE, RMSE, RMSLE and MAXE of the models, namely PPR78,
PPR78+Adj \citep{AlGhafri2014},  PPR78+InG and PPR78+Opt (present work), presented in Table \ref{table 2:Maen Of Metrics}. For each approach,
the metric values obtained, for Live oil 1 and Live oil 2, and also for various temperatures 323.15 K, 373.15 K and 423.15 K are averaged and presented in
Table \ref{table 2:Maen Of Metrics}.
Using the same five evaluation measures as before (described in section \ref{Metrics}) are utilized for each approach, with the best results shown in bold. As can be seen, according to all metrics reported in Table \ref{table 2:Maen Of Metrics}, our suggested model PPR78+Opt performs noticeably better than predictive methods PPR78, PPR78+Adj  [\textbf{1}] and even PPR78+InG introduced in present work. Particularly,  mean of the cost function of PPR78+Opt with respect to MSE, RMSE and RMSLE metrics is  0.911, 0.912 and 0.037, respectively,  which are significantly lower than the results of other predictive models. It is interesting to note that even PPR78+Adj, which uses an adjustment strategy to lower the cost function, has an MSE of 9.004, which is considerably worse than our proposed method. This demonstrates how well our suggested optimization model performs in comparison to other prediction models. In order to visualize this superiority, the result of Table \ref{table 2:Maen Of Metrics} is presented by a spider-web diagram (see Figure \ref{fig:spiderweb}), which is suitable to compare the effectiveness of various methods. As shown in Figure \ref{fig:spiderweb}, and in agreement with the results from Tables \ref{table1 :costFunctions} and \ref{table 2:Maen Of Metrics} in this section, PPR78+Opt is the best method, with the lowest cost function in relation to all the metrics. In order to scale the results so that they may fit in a single spider-web figure, the mean of cost function values of all the methods with regard to the metrics MSE and RMSLE are divided and multiplied by 10, respectively.

As discussed earlier, a slight agreement was seen at the experimental bubble points by adjusting the single parameter, binary interaction of $CO_2$-$CH_4$ at each temperature, but quite significant discrepancies remained unchanged at the dew points,  (see red dashed and solid lines in Figure \ref{fig:6fig}).  The values of this parameter as predicted by PPR78, adjusted by PPR78+Adj, and optimized by PPR78+Opt are reported in Table \ref{Table3}. Furthermore, Figure \ref{kijComp} represents information from Table \ref{Table3} in an illustrative manner.  As we can see, the binary interaction adjusted  [\textbf{1}] deviates significantly from the values predicted by PPR78.  
In contrast, the binary interaction suggested by our optimization method is reasonably close to the values suggested by the predictive method (PPR78), and results in reduction of MSE cost function. Additionally, we can see the agreement between the binary interaction we suggested and the outcomes of the predictive method by looking at the binary interaction computed for $CO_2$ with other components. Thus, the PPR78+Opt is the best method to find the optimal value for binary interaction between $CO_2$ and $CH_4$ and results in reduce the cost function, according to the results already shown and analyzed in this Section.

\begin{table}[h]
	\begin{minipage}[c]{.44\textwidth }%
		\footnotesize\centering
		\renewcommand{\arraystretch}{1}
		\captionof{table}{\small Binary interaction parameters for $CO_2$–$CH_4$, predicted by PPR78, adjusted by PPR78+Adj and optimized by PPR78+Opt;}
		
		\begin{tabular}{|c|c c|c|}
			\hline
			\multirow{2}{*}{Temp [K]}  & \multicolumn{3}{c|}{$k_{CO_2-CH_4}$}\\ 
			\cline{2-4}
			&\multicolumn{1}{l|}{PPR78}  & \multicolumn{1}{l|}{ PPR78+Adj} & PPR78+Opt \\
			\hline
			323.15 &\multicolumn{1}{l|}{0.119} &  -0.07 & 0.116   \\
			\hline
			373.15 & \multicolumn{1}{l|}{0.13} & -0.2 & 0.105  \\
			\hline
			423.15 & \multicolumn{1}{l|}{0.142} & -0.25 & 0.107 \\  
			\hline
			
		\end{tabular}
		\label{Table3}  
		

		\hrule height 0pt
	\end{minipage}%
	\begin{minipage}[c]{.7\textwidth}
		\centering
		\includegraphics[scale=0.42, trim={0 0 0 0}]{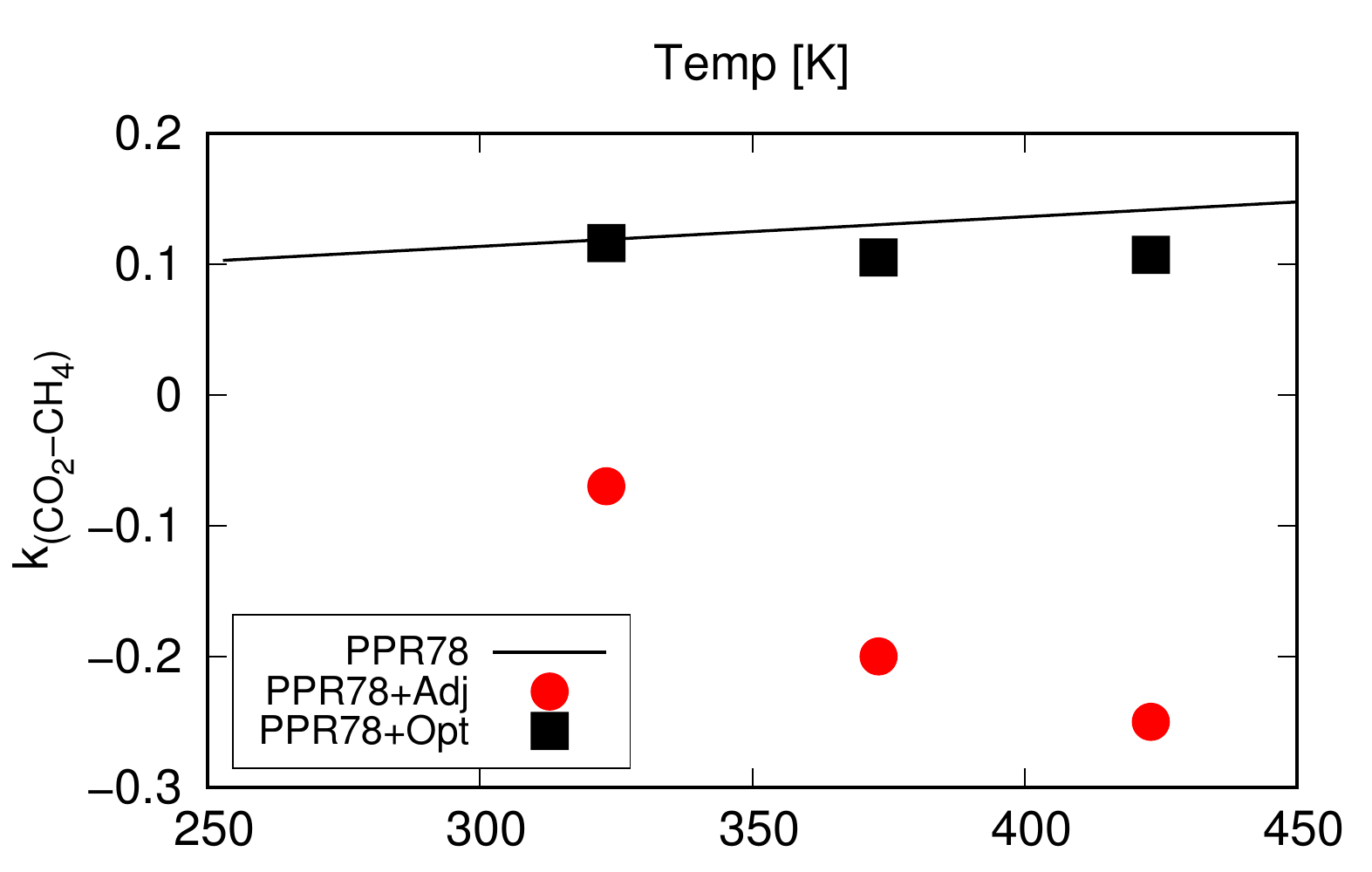}
		
		\captionof{figure}{\small  Binary interaction \\ parameters reported in Table \ref{Table3}}
		\label{kijComp}
	\end{minipage}
\end{table}

\subsection{Concluding remarks and future works}

This work employs the grid search technique to optimize the predictions of the phase
behavior of systems containing $CO_2$ and multi-component hydrocarbons at
reservoir temperatures and pressures. Using some physical intuition, it was realized before, that
adjusting the binary interaction parameter between $CO_2$ and $CH_4$ improves the
results. However, by adjusting this single parameter at each temperature, a slight agreement
was observed at the experimental bubble points, while quite big discrepancies remained at the
dew points.

We improve the results (reducing the MSE) by employing a reasonable initial guess strategy and then find the optimal solution due to single evaluation task (either positive or negative $k_{CO_2-CH_4}$), which is the cost function minimizer (see Tables \ref{table1 :costFunctions}, \ref{table 2:Maen Of Metrics} and Figures \ref{fig:6fig}, \ref{fig:spiderweb}), we observe that
our proposed optimization approach (PPR78+Opt) accurately predicts the experimental results in either the dew or bubble points. This clarify that PPR78+Opt provides the smallest cost function when compared to other predictive methods. It should be noted that this results is due to grid search strategy is computationally heavy. 

\begin{figure}[h]   
	\centering
	\includegraphics[ scale=0.18, trim={5cm 0cm 3cm 0cm}]{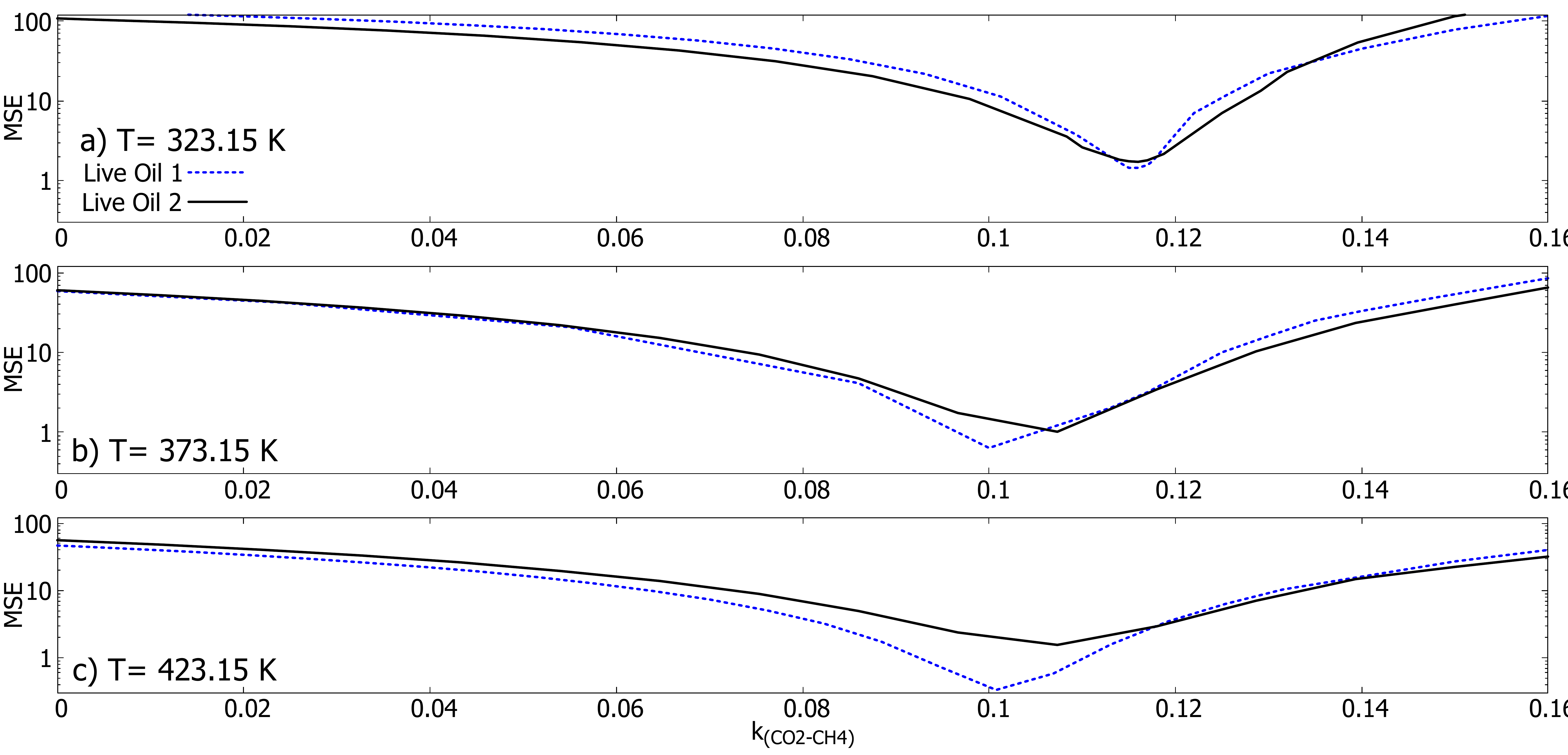}
	\captionof{figure}{\small MSE cost functions obtained by PPR78+Opt }
	\label{fig:3figs}
\end{figure}

The next step that can be done in the future is based on what we discovered from the cost function (MSE) of the PPR78+Opt model (see Figure. \ref{fig:3figs}); we discovered that MSE behaved in a way that can be classified as a "One Dimensional Optimization Problem".  To put it another way, the cost function is a Unimodal function with a convex epigraph and a unique global minimum. 
Thus, we can use a straightforward optimization strategies (Fibonacci Search Method or Golden Section Search Method) to reduce CPU time by finding the optimum $k_{CO_2-CH_4}$  rather than grid search to obtain the smallest MSE. This leads to reducing the cost of computation, with out affecting the cost function results. This improvement can apply without any complication and improvement in CPU time is ensured by the superiority of one dimensional search methods over the grid search.

As a further stage, it is reasonable to attempt to use an optimization strategy that improves decision making, such as Bayesian Optimization (BO), to address the issue of inaccurate binary interaction prediction; it is possible to adjust other parameters (around 200 parameters) blindly. Thus, the BO method decides which parameter (binary interaction between $CO_2$
and other components) needs to be optimized. The target is to supply the solver with a few
experimental data, and the method automatically adjusts the parameters.

\section*{Acknowledgments}
The authors acknowledge the support awarded by GALP/Petrogal Brazil and ANP–Agência Nacional de Petróleo, Gás Natural e Biocombustível.

\section*{References}
\begin{enumerate}[label={[{\textbf{\arabic*}}]}]
	\item
	Saif Z. AlGhafri, Geoffrey C. Maitland, and J.P. Martin Trusler. Experimental
	and modeling study of the phase behavior of synthetic crude oil+co2. Fluid
	Phase Equilibria, 365:20-40, 3 2014.
	\item 
	João A P. Coutinho, Georgios M. Kontogeorgis, and Erling H. Stenby. Binary
	interaction parameters for nonpolar systems with cubic equations of state: a
	theoretical approach 1. co2/hydrocarbons using srk equation of state. Fluid
	Phase Equilibria, 102(1):31-60, 1994.
	\item 
	Jean-Noël. Jaubert and Fabrice. Mutelet. Vle predictions with the peng-Robinson
	equation of state and temperature dependent kij calculated through a group
	contribution method. Fluid Phase Equilibria, 224:285-304, 10 2004.
	\item 
	Jean-Noël. Jaubert, Romain. Privat, and Fabrice. Mutelet. Predicting the phase
	equilibria of synthetic petroleum fluids with the PPR78 approach. AIChE Journal,
	56:3225-3235, 12 2010.
	\item 
	Jean-Noël. Jaubert, Stéphane. Vitu, Fabrice. Mutelet, and Jean-Pierre. Corriou.
	Extension of the PPR78 model (predictive 1978, peng-Robinson EOS with temperature
	dependent kij calculated through a group contribution method) to systems
	containing aromatic compounds. Fluid Phase Equilibria, 237:193-211, 10 2005.
	\item 
	Michael L. Michelsen. The isothermal flash problem. part i. stability. Fluid
	Phase Equilibria, 9:1-19, 12 1982.
	\item 
	Michael L. Michelsen and Jørgen M. Mollerup. Thermodynamic models : fundamentals $\&$ computational aspects. Tie-Line Publications, 2007.
	\item 
	Ding-Yu. Peng and Donald B. Robinson. A new two-constant equation of state.
	Industrial $\&$ Engineering Chemistry Fundamentals, 15:59-64, 2 1976.
	\item 
	Iuri S V. Segtovich, Amaro G. Barreto, and Frederico W. Tavares. Simultaneous
	multiphase flash and stability analysis calculations including hydrates. Fluid
	Phase Equilibria, 413:196-208, 4 2016.
 \end{enumerate}

\bibliographystyle{plain}
\end{document}